\documentstyle[12pt]{article}
\begin{document}

\title{On continuous Polish group actions and equivalence relations}

\author{Nikolaos E. Sofronidis\footnote{$A \Sigma MA:$ 130/2543/94}}

\date{\footnotesize Department of Economics, University of Ioannina, Ioannina 45110, Greece.
(nsofron@otenet.gr, nsofron@cc.uoi.gr)}

\maketitle

\begin{abstract}
Let $X = \left\{ P \in [0,1]^{\bf N} : \left( \forall \nu \in {\bf
N} \right) \left( P \left( \{ \nu \} \right) > 0 \right) \wedge
\sum\limits_{ \nu = 0 }^{ \infty } P \left( \{ \nu \} \right) = 1
\right\} $ be the Polish space of probability measures on ${\bf
N}$, each of which assigns positive probability to every
elementary event, while for any $P \in X$, let ${\Gamma }_{P} =
\left\{ \xi \in L^{1}( {\bf N} , P ) : \left( \forall \nu \in {\bf
N} \right) \left( \xi ( \nu ) > 0 \right) \wedge \sum\limits_{ \nu
= 0 }^{ \infty } \xi ( \nu ) P \left( \{ \nu \} \right) = 1
\right\} $ and let ${\Phi }_{P} : {\Gamma }_{P} \ni \xi \mapsto
{\Phi }_{P}( \xi ) \in X$ be defined by the relation $\left( {\Phi
}_{P}( \xi ) \right) \left( \{ \nu \} \right) = \xi ( \nu ) P
\left( \{ \nu \} \right) $, whenever $\nu \in {\bf N}$. If we
consider the equivalence relation $E = \left\{ (P,Q) \in X^{2} :
\left( \exists \xi \in {\Gamma }_{P} \right) \left( Q = {\Phi
}_{P}( \xi ) \right) \right\} $, the Polish space ${\bf P} =
\left\{ {\bf x} \in {\ell }^{1} \left( {\bf R} \right) : \left(
\forall n \in {\bf N} \right) \left( {\bf x}(n) > 0 \right)
\right\} $ and the commutative Polish group ${\bf G} = \left\{
{\bf g} \in (0 , \infty )^{\bf N} : \lim\limits_{n \rightarrow
\infty }{\bf g}(n) = 1 \right\} $, while we set $\left( {\bf g}
\cdot {\bf x} \right) (n) = {\bf g}(n){\bf x}(n)$, whenever ${\bf
g} \in {\bf G}$, ${\bf x} \in {\bf P}$ and $n \in {\bf N}$, then
$E$ is definable and it admits a strong approximation by the
turbulent Polish group action of ${\bf G}$ on ${\bf P}$. In
addition, if we consider the Polish space ${\ell }^{1} \left( {\bf
C}^{*} \right) = \left\{ {\bf x} \in {\ell }^{1} \left( {\bf C}
\right) : \left( \forall n \in {\bf N} \right) \left( {\bf x}(n)
\neq 0 \right) \right\} $ and we set ${\bf H} = \left\{ {\bf h}
\in \left( {\bf C}^{*} \right) ^{\bf N} : \lim\limits_{n
\rightarrow \infty }{\bf h}(n) = 1 \right\} $, while $\left( {\bf
h} \cdot {\bf x} \right) (n) = {\bf h}(n){\bf x}(n)$, whenever
${\bf h} \in {\bf H}$, ${\bf x} \in {\ell }^{1} \left( {\bf C}^{*}
\right) $ and $n \in {\bf N}$, then ${\bf H}$ is a commutative
Polish group under pointwise multiplication and ${\bf H} \times
{\ell }^{1}\left( {\bf C}^{*} \right) \ni \left( {\bf h} , {\bf x}
\right) \mapsto {\bf h} \cdot {\bf x} \in {\ell }^{1}\left( {\bf
C}^{*} \right) $ constitutes a continuous Polish group action each
orbit of which is dense and meager, while ${\bf G}$ on ${\bf P}$
is a subaction of ${\bf H}$ on ${\ell }^{1} \left( {\bf C}^{*}
\right) $. In addition, if $\lambda $ is the Lebesgue measure in
the real line and we consider the Polish space $L^{1}_{++} \left(
(0, \infty ) , \lambda \right) = \left\{ f \in L^{1} \left( (0,
\infty ) , \lambda \right) : f > 0, \ \lambda -a.e. \right\} $,
while we consider the set ${\bf F} = \left\{ f \in C \left( (0,
\infty ) , (0, \infty ) \right) : \lim\limits_{x \rightarrow
0}f(x) = \lim\limits_{x \rightarrow \infty }f(x) = 1 \right\} $
and the operation $\left( f \cdot g \right) (x) = f(x)g(x)$,
whenever $f \in {\bf F}$, $g \in L^{1}_{++} \left( (0, \infty ) ,
\lambda \right) $ and $x \in (0, \infty )$, then ${\bf F}$
constitutes a commutative Polish group under pointwise
multiplication and ${\bf F} \times L^{1}_{++} \left( (0, \infty )
, \lambda \right) \ni (f,g) \mapsto f \cdot g \in L^{1}_{++}
\left( (0, \infty ) , \lambda \right) $ constitutes a continuous
Polish group action, which is not an extension of the turbulent
Polish group action of the group ${\bf G}^{*} = \left\{ {\bf g}
\in (0 , \infty )^{{\bf N}^{*}} : \lim\limits_{n \rightarrow
\infty }{\bf g}(n) = 1 \right\} $, which is essentially ${\bf G}$,
on the space ${\bf P}^{*} = \left\{ {\bf x} \in (0, \infty )^{{\bf
N}^{*}} : \sum\limits_{n=1}^{\infty }{\bf x}(n) < \infty \right\}
$, which is essentially ${\bf P}$, even though ${\bf G}^{*}$,
${\bf P}^{*}$ are Polish subspaces of ${\bf F}$, $L^{1}_{++}
\left( ( 0 , \infty ) , \lambda \right) $ respectively.
\end{abstract}

\section*{\footnotesize{{\bf Mathematics Subject Classification:} 03E15, 22A99, 40A99.}}

\section{Introduction}

The following is a modification of the one given in [2].
\\ \rm \\
{\bf 1.1. Definition.} If $X$ is any Polish space and $E$ is any
definable equiva- lence relation on $X$, then $E$ admits a strong
approximation by a Polish group action, when the following
conditions are satisfied:
\begin{enumerate}
\item[(i)]
For any $x \in X$, there exists a Polish space ${\Gamma }_{x}$ and
a continuous mapping ${\phi }_{x} : {\Gamma }_{x} \rightarrow X$
such that $x \mapsto {\Gamma }_{x}$ and $x \mapsto {\phi }_{x}$
are definable, while ${\phi }_{x} \left[ {\Gamma }_{x} \right] =
[x]_{E}$.
\item[(ii)]
There exists a Polish group $G$ acting continuously on a Polish
space $P$ with the property that $X \subseteq P$ and for any $x
\in X$, there exists a Polish space ${\Delta }_{x}$ such that
${\Gamma }_{x}$ is a Polish subspace of ${\Delta }_{x}$ and $x
\mapsto {\Delta }_{x}$ is definable, while ${\Gamma }_{x}
\subseteq {\overline{ {\Delta }_{x} \cap G}}$ and ${\phi }_{x} (g)
= g \cdot x$, whenever $g \in {\Gamma }_{x} \cap G$.
\item[(iii)]
For any $x \in X$ and for any $\gamma \in {\Gamma }_{x}$, there
exists a homeomorphism  ${\psi }_{x , \gamma } : {\Gamma }_{x}
\rightarrow {\Gamma }_{ {\phi }_{x}( \gamma ) }$ with the property
that $( x, \gamma ) \mapsto {\psi }_{x , \gamma }$ is definable
and ${\phi }_{x}( \delta ) = {\phi }_{ {\phi }_{x}( \gamma ) }
\left( {\psi }_{x, \gamma }( \delta ) \right) $, whenever $\delta
\in {\Gamma }_{x}$.
\end{enumerate}

\noindent The following is introduced in [1].
\\ \rm \\
{\bf 1.2. Definition.} Let $G$ be any Polish group with countable
base ${\mathcal{B}}_{G}$ acting continuously on a Polish space $X$
with countable base ${\mathcal{B}}_{X}$ and let $$x E^{X}_{G} y
\iff \left( \exists g \in G \right) \left( g \cdot x = y \right)
,$$ whenever $x$, $y$ are in $X$. For any open neighborhood $U \in
{\mathcal{B}}_{X}$ of $x$ and for any symmetric open neighborhood
$V \in {\mathcal{B}}_{G}$ of $1^{G}$, the $(U,V)$-local orbit
$O(x,U,V)$ of $x$ in $X$ is defined, as follows:
\begin{enumerate}
\item[]
$y \in O(x,U,V)$ if there exist $g_{0}$, ..., $g_{k}$ in $V$,
where $k \in {\bf N}$, such that if $x_{0} = x$ and $x_{i+1} =
g_{i} \cdot x_{i}$ for every $i \in \{ 0, ..., k \} $, then all
the $x_{i}$ are in $U$ and $x_{k+1} = y$.
\end{enumerate}
The action of $G$ on $X$ is said to be turbulent at the point $x$,
in symbols $x \in T^{X}_{G}$, if for any such $U$ and $V$, there
exists an open neighborhood $U' \in {\mathcal{B}}_{X}$ of $x$ such
that $U ' \subseteq U$ and $O(x,U,V)$ is dense in $U'$.
\\ \rm \\
{\bf 1.3. Definition.} Let $$X = \left\{ P \in [0,1]^{\bf N} :
\left( \forall \nu \in {\bf N} \right) \left( P \left( \{ \nu \}
\right) > 0 \right) \wedge \sum\limits_{ \nu = 0 }^{ \infty } P
\left( \{ \nu \} \right) = 1 \right\} $$ be the Polish space of
probability measures on ${\bf N}$, each of which assigns positive
probability to every elementary event. (It is a $G_{\delta }$
subset of the compact Polish space $[0,1]^{\bf N}$ equipped with
the product topology.) For any $P \in X$, let $${\Gamma }_{P} =
\left\{ \xi \in L^{1}( {\bf N} , P ) : \left( \forall \nu \in {\bf
N} \right) \left( \xi ( \nu ) > 0 \right) \wedge \sum\limits_{ \nu
= 0 }^{ \infty } \xi ( \nu ) P \left( \{ \nu \} \right) = 1
\right\} $$ be the Polish space of positive $L^{1}$ random
variables on the probability measure space $\left( {\bf N} ,
{\mathcal{P}}( {\bf N} ) , P \right) $ whose expectation is equal
to $1$ (it is a $G_{\delta }$ subset of the $C2$ Banach space
$L^{1}( {\bf N} , P )$) and let ${\Phi }_{P} : {\Gamma }_{P} \ni
\xi \mapsto {\Phi }_{P}( \xi ) \in X$ be defined by the relation
$\left( {\Phi }_{P}( \xi ) \right) \left( \{ \nu \} \right) = \xi
( \nu ) P \left( \{ \nu \} \right) $, whenever $\nu \in {\bf N}$.
We set $$E = \left\{ (P,Q) \in X^{2} : \left( \exists \xi \in
{\Gamma }_{P} \right) \left( Q = {\Phi }_{P}( \xi ) \right)
\right\} .$$
\\
The following is proved in [3].
\\ \rm \\
{\bf 1.4. Theorem.} If ${\bf P} = \left\{ {\bf x} \in {\ell }^{1}
\left( {\bf R} \right) : \left( \forall n \in {\bf N} \right)
\left( {\bf x}(n)
> 0 \right) \right\} $ and $${\bf G} = \left\{ {\bf g} \in (0 ,
\infty )^{\bf N} : \lim\limits_{n \rightarrow \infty }{\bf g}(n) =
1 \right\} ,$$ while $\left( {\bf g} \cdot {\bf x} \right) (n) =
{\bf g}(n){\bf x}(n)$, whenever ${\bf g} \in {\bf G}$, ${\bf x}
\in {\bf P}$ and $n \in {\bf N}$, then the following hold true:
\begin{enumerate}
\item[(i)]
${\bf G}$ constitutes a commutative Polish group under pointwise
multiplication.
\item[(ii)]
${\bf G} \times {\bf P} \ni \left( {\bf g} , {\bf x} \right)
\mapsto {\bf g} \cdot {\bf x} \in {\bf P}$ constitutes a
continuous Polish group action.
\item[(iii)]
The action of ${\bf G}$ on ${\bf P}$ is turbulent.
\end{enumerate}

\noindent So our first purpose in this article is to prove the
following.
\\ \rm \\
{\bf 1.5. Theorem.} $E$ is definable and it admits a strong
approximation by the turbulent Polish group action of ${\bf G}$ on
${\bf P}$.
\\ \rm \\
{\bf 1.6. Definition.} If ${\bf C}^{*} = {\bf C} \setminus \{ 0 \}
$, then we set $${\ell }^{1} \left( {\bf C}^{*} \right) = \left\{
{\bf x} \in {\ell }^{1} \left( {\bf C} \right) : \left( \forall n
\in {\bf N} \right) \left( {\bf x}(n) \neq 0 \right) \right\} ,$$
which is a $G_{\delta }$ subset of the $C2$ Banach space $${\ell
}^{1} \left( {\bf C} \right) = \left\{ {\bf x} \in {\bf C}^{\bf N}
: \sum\limits_{n=0}^{ \infty } \left\vert {\bf x}(n) \right\vert <
\infty \right\} $$ and consequently it constitutes a Polish space.
\\ \rm \\
So our second purpose in this article is to prove the following.
\\ \rm \\
{\bf 1.7. Theorem.} If $${\bf H} = \left\{ {\bf h} \in \left( {\bf
C}^{*} \right) ^{\bf N} : \lim\limits_{n \rightarrow \infty }{\bf
h}(n) = 1 \right\} ,$$ while $\left( {\bf h} \cdot {\bf x} \right)
(n) = {\bf h}(n){\bf x}(n)$, whenever ${\bf h} \in {\bf H}$, ${\bf
x} \in {\ell }^{1} \left( {\bf C}^{*} \right) $ and $n \in {\bf
N}$, then the following hold true:
\begin{enumerate}
\item[(i)]
${\bf H}$ constitutes a commutative Polish group under pointwise
multiplication and ${\bf G}$ is a Polish subgroup of ${\bf H}$.
\item[(ii)]
${\bf H} \times {\ell }^{1}\left( {\bf C}^{*} \right) \ni \left(
{\bf h} , {\bf x} \right) \mapsto {\bf h} \cdot {\bf x} \in {\ell
}^{1}\left( {\bf C}^{*} \right) $ constitutes a continuous Polish
group action and ${\bf G} \times {\bf P} \ni \left( {\bf g} , {\bf
x} \right) \mapsto {\bf g} \cdot {\bf x} \in {\bf P}$ is a
subaction.
\item[(iii)]
For any ${\bf x} \in {\ell }^{1} \left( {\bf C}^{*} \right) $,
${\bf H} \cdot {\bf x}$ is dense in ${\ell }^{1} \left( {\bf
C}^{*} \right) $.
\item[(iv)]
For any ${\bf x} \in {\ell }^{1} \left( {\bf C}^{*} \right) $,
${\bf H} \cdot {\bf x}$ is meager in ${\ell }^{1} \left( {\bf
C}^{*} \right) $.
\end{enumerate}

\noindent The following Polish space is due to Tom Wolff and is
introduced in [2]. It is a $G_{\delta }$ subset of a closed subset
of a Polish space. See (i) of Proposition 5.6 on page 1470 of [2].
\\ \rm \\
{\bf 1.8. Definition.} If $\lambda $ is the Lebesgue measure in
the real line, then we set $L^{1}_{++} \left( (0, \infty ) ,
\lambda \right) = \left\{ f \in L^{1} \left( (0, \infty ) ,
\lambda \right) : f > 0, \ \lambda -a.e. \right\} $.
\\ \rm \\
So, following [3], our third purpose in this article is to prove
the following.
\\ \rm \\
{\bf 1.9. Theorem.} If ${\bf F} = \left\{ f \in C \left( (0,
\infty ) , (0, \infty ) \right) : \lim\limits_{x \rightarrow
0}f(x) = \lim\limits_{x \rightarrow \infty }f(x) = 1 \right\} $
and $\left( f \cdot g \right) (x) = f(x)g(x)$, whenever $f \in
{\bf F}$, $g \in L^{1}_{++} \left( (0, \infty ) , \lambda \right)
$ and $x \in (0, \infty )$, then the following hold true:
\begin{enumerate}
\item[(i)]
${\bf F}$ constitutes a commutative Polish group under pointwise
multiplication.
\item[(ii)]
${\bf F} \times L^{1}_{++} \left( (0, \infty ) , \lambda \right)
\ni (f,g) \mapsto f \cdot g \in L^{1}_{++} \left( (0, \infty ) ,
\lambda \right) $ constitutes a continuous Polish group action.
\end{enumerate}

\noindent An immediate consequence of [3] is the following.
\\ \rm \\
{\bf 1.10. Theorem.} If ${\bf P}^{*} = \left\{ {\bf x} \in (0,
\infty )^{{\bf N}^{*}} : \sum\limits_{n=1}^{\infty }{\bf x}(n) <
\infty \right\} $ and $${\bf G}^{*} = \left\{ {\bf g} \in (0 ,
\infty )^{{\bf N}^{*}} : \lim\limits_{n \rightarrow \infty }{\bf
g}(n) = 1 \right\} ,$$ while $\left( {\bf g} \cdot {\bf x} \right)
(n) = {\bf g}(n){\bf x}(n)$, whenever ${\bf g} \in {\bf G}^{*}$,
${\bf x} \in {\bf P}^{*}$ and $n \in {\bf N}^{*}$, then the
following hold true:
\begin{enumerate}
\item[(i)]
${\bf G}^{*}$ constitutes a commutative Polish group under
pointwise multiplication.
\item[(ii)]
${\bf G}^{*} \times {\bf P}^{*} \ni \left( {\bf g} , {\bf x}
\right) \mapsto {\bf g} \cdot {\bf x} \in {\bf P}^{*}$ constitutes
a continuous Polish group action.
\item[(iii)]
The action of ${\bf G}^{*}$ on ${\bf P}^{*}$ is turbulent.
\end{enumerate}

\noindent {\bf 1.11. Definition.} If ${\bf g} \in {\bf G}^{*}$,
then we set
\[
f_{\bf g}(x) = \left\{
\begin{array}{lllll}
\left( {\bf g}(n+1) - {\bf g}(n) \right) (x-n) + {\bf g}(n) &
\mbox{if $x \in [n,n+1]$ and $n \in {\bf N}^{*}$}
\\ \\
2 \left( {\bf g}(1) - 1 \right) (x-1) + {\bf g}(1) & \mbox{if $x
\in \left[ \frac{1}{2} , 1 \right] $}
\\ \\
1 & \mbox{if $x \in \left( \left. 0 , \frac{1}{2} \right] \right.
$}
\end{array}
\right.
\]
and it is not difficult to verify that $f_{\bf g} \in C \left( ( 0
, \infty ) , ( 0 , \infty ) \right) $ due to the fact that ${\bf
g} \in ( 0 , \infty )^{{\bf N}^{*}}$ and $( 0 , \infty )^{2}$ is
convex, while obviously $\lim\limits_{x \rightarrow 0} f_{\bf
g}(x) = 1$ and $\lim\limits_{x \rightarrow \infty } f_{\bf
g}(x)=1$, since $\lim\limits_{n \rightarrow \infty } {\bf g}(n) =
1$ and $( 0 , \infty )^{2}$ is convex, i.e., $f_{\bf g} \in {\bf
F}$.
\\ \rm \\
So, our fourth purpose in this article is to prove the following.
\\ \rm \\
{\bf 1.12. Theorem.} ${\bf G}^{*} \ni {\bf g} \mapsto f_{\bf g}
\in {\bf F}$ is a Polish space continuous injection, but not a
Polish group continuous injection.
\\ \rm \\
{\bf 1.13. Remark.} If ${\bf x} \in {\bf P}^{*}$, then we set
\[
g_{\bf x}(x) = \left\{
\begin{array}{lll}
{\bf x}(n) & \mbox{if $x \in [n,n+1)$ and $n \in {\bf N}^{*}$}
\\ \\
1 & \mbox{if $x \in (0,1)$}
\end{array}
\right.
\]
and it is not difficult to verify that $g_{\bf x} \in L^{1}_{++}
\left( ( 0 , \infty ) , \lambda \right) $ and $\left\Vert g_{\bf
x} \right\Vert _{1} = 1 + \Vert {\bf x} \Vert _{1}$, so ${\bf
P}^{*} \ni {\bf x} \mapsto g_{\bf x} \in L^{1}_{++} \left( ( 0 ,
\infty ) , \lambda \right) $ is a Polish space continuous
injection.

\section{The proof of 1.5}

\subsection{The proof of (i) in 1.1}

If ${\xi }^{k} \rightarrow \xi $ in ${\Gamma }_{P}$ as $k
\rightarrow \infty $, then for any $\nu \in {\bf N}$, we have that
\[
\begin{array}{llll}
\left\vert \left( {\Phi }_{P} \left( {\xi }^{k} \right) \right)
\left( \{ \nu \} \right) - \left( {\Phi }_{P} \left( \xi \right)
\right) \left( \{ \nu \} \right) \right\vert & = \ \left\vert {\xi
}^{k} \left( \{ \nu \} \right) P \left( \{ \nu \} \right) - \xi (
\nu ) P \left( \{ \nu \} \right) \right\vert
\\ \\
 & = \ \left\vert {\xi }^{k}( \nu ) - \xi ( \nu ) \right\vert P \left( \{ \nu \} \right)
\\ \\
 & \leq \ \sum\limits_{n=0}^{ \infty } \left\vert {\xi }^{k}(n) - \xi (n) \right\vert P \left( \{ n \} \right)
\\ \\
 & = \ \left\Vert {\xi }^{k} - \xi \right\Vert _{1} \rightarrow 0
\end{array}
\]
as $k \rightarrow \infty $, so for any $\nu \in {\bf N}$, we have
that $$\lim\limits_{k \rightarrow \infty } \left( {\Phi }_{P}
\left( {\xi }^{k} \right) \right) \left( \{ \nu \} \right) =
\left( {\Phi }_{P} \left( \xi \right) \right) \left( \{ \nu \}
\right) $$ and consequently ${\Phi }_{P} \left( {\xi }^{k} \right)
\rightarrow {\Phi }_{P} \left( \xi \right) $ in $X$ as $k
\rightarrow \infty $. Moreover, if $(P,Q) \in E$ and $\xi \in
{\Gamma }_{P}$ is such that $Q = {\Phi }_{P} \left( \xi \right) $,
then $\frac{1}{ \xi } \in {\Gamma }_{Q}$ and $P = {\Phi }_{Q}
\left( \frac{1}{ \xi } \right) $. Indeed, $\frac{1}{ \xi ( \nu ) }
> 0$ for every $\nu \in {\bf N}$ and $\sum\limits_{ \nu = 0 }^{
\infty } \frac{1}{ \xi ( \nu ) } Q \left( \{ \nu \} \right) =
\sum\limits_{ \nu = 0 }^{ \infty } \frac{1}{ \xi ( \nu ) } {\Phi
}_{P} \left( \{ \nu \} \right) = \sum\limits_{ \nu = 0 }^{ \infty
} \frac{1}{ \xi ( \nu ) } \xi ( \nu ) P \left( \{ \nu \} \right) =
1$, hence $\frac{1}{ \xi } \in {\Gamma }_{Q}$, while $\left( {\Phi
}_{Q} \left( \frac{1}{ \xi } \right) \right) \left( \{ \nu \}
\right) = \frac{1}{ \xi ( \nu ) } Q \left( \{ \nu \} \right) =
\frac{1}{ \xi ( \nu ) } {\Phi }_{P} \left( \{ \nu \} \right) =
\frac{1}{ \xi ( \nu ) } \xi ( \nu ) P \left( \{ \nu \} \right) = P
\left( \{ \nu \} \right) $ for every $\nu \in {\bf N}$. \hfill
$\bigtriangleup $

\subsection{The proof of (ii) in 1.1}

We shall first prove that ${\Gamma }_{P} \subseteq {\overline{
L^{1}( {\bf N} , P ) \cap {\bf G}}}$. Indeed, if $\xi \in {\Gamma
}_{P}$ and $\epsilon > 0$, then there exists $N \in {\bf N}$ such
that $\sum\limits_{ \nu > N } P \left( \{ \nu \} \right) < \frac{
\epsilon }{2}$ and $\sum\limits_{ \nu > N } \xi ( \nu ) P \left(
\{ \nu \} \right) < \frac{ \epsilon }{2}$, so if ${\bf g} = \left(
\xi (0) , ..., \xi (N), 1, 1, ... \right) $, then ${\bf g} \in
L^{1}( {\bf N} , P ) \cap {\bf G}$ and $\left\Vert \xi - {\bf g}
\right\Vert _{1} = \sum\limits_{ \nu > N } \left\vert \xi ( \nu )
- 1 \right\vert P \left( \{ \nu \} \right) \leq \sum\limits_{ \nu
> N } \xi ( \nu ) P \left( \{ \nu \} \right) + \sum\limits_{ \nu >
N } P \left( \{ \nu \} \right) < \frac{ \epsilon }{2} + \frac{
\epsilon }{2} = \epsilon $. Moreover, for any ${\bf g} \in {\Gamma
}_{P} \cap {\bf G}$, we have that ${\Phi }_{P}\left( {\bf g}
\right) = {\bf g} \cdot P$, since $\left( {\Phi }_{P} \left( {\bf
g} \right) \right) \left( \{ \nu \} \right) = {\bf g}( \nu ) P
\left( \{ \nu \} \right) = \left( {\bf g} \cdot P \right) ( \nu )$
for every $\nu \in {\bf N}$. \hfill $\bigtriangleup $

\subsection{The proof of (iii) in 1.1}

If $P \in X$ and $\xi \in {\Gamma }_{P}$, then we define ${\Psi
}_{P , \xi } : {\Gamma }_{P} \ni \zeta \mapsto {\Psi }_{P, \xi }
\left( \zeta \right) \in {\Gamma }_{ {\Phi }_{P} \left( \xi
\right) }$ by the relation $\left( {\Psi }_{P , \xi } \left( \zeta
\right) \right) ( \nu ) = \frac{ \zeta ( \nu ) }{ \xi ( \nu ) }$,
whenever $\nu \in {\bf N}$. It is not difficult to see that
$\sum\limits_{ \nu = 0 }^{ \infty } \left( {\Psi }_{P , \xi }
\left( \zeta \right) \right) ( \nu ) \left( {\Phi }_{P} \left( \xi
\right) \right) \left( \{ \nu \} \right) = \sum\limits_{ \nu = 0
}^{ \infty } \frac{ \zeta ( \nu ) }{ \xi ( \nu ) } \xi ( \nu ) P
\left( \{ \nu \} \right) = \sum\limits_{ \nu = 0 }^{ \infty }
\zeta ( \nu ) P \left( \{ \nu \} \right) = 1$ and $\frac{ \zeta (
\nu ) }{ \xi ( \nu ) } > 0$ for every $\nu \in {\bf N}$, so ${\Psi
}_{P , \xi }$ is well-defined. In addition, if ${\zeta }_{1}$,
${\zeta }_{2}$ are any elements of ${\Gamma }_{P}$, then
\begin{enumerate}
\item[ ]
$\left\Vert {\Psi }_{P , \xi } \left( {\zeta }_{1} \right) - {\Psi
}_{P , \xi } \left( {\zeta }_{2} \right) \right\Vert _{1}$
\item[ ]
$= \sum\limits_{ \nu = 0 }^{ \infty } \left\vert \left( {\Psi }_{P
, \xi } \left( {\zeta }_{1} \right) \right) ( \nu ) - \left( {\Psi
}_{P , \xi } \left( {\zeta }_{2} \right) \right) ( \nu )
\right\vert \left( {\Phi }_{P} \left( \xi \right) \right) \left(
\{ \nu \} \right) $
\item[ ]
$= \sum\limits_{ \nu = 0 }^{ \infty } \left\vert \frac{ {\zeta
}_{1}( \nu ) }{ \xi ( \nu ) } - \frac{ {\zeta }_{2}( \nu ) }{ \xi
( \nu ) } \right\vert \xi ( \nu ) P \left( \{ \nu \} \right) $
\item[ ]
$= \sum\limits_{ \nu = 0 }^{ \infty } \left\vert {\zeta }_{1}( \nu
) - {\zeta }_{2}( \nu ) \right\vert P \left( \{ \nu \} \right) $
\item[ ]
$= \left\Vert {\zeta }_{1} - {\zeta }_{2} \right\Vert _{1}$
\end{enumerate}
and ${\Psi }_{P , \xi }^{-1} = {\Psi }_{P , \frac{1}{ \xi } }$,
since $\left( {\Psi }_{P , \frac{1}{ \xi } } \left( {\Psi }_{P ,
\xi } \left( \zeta \right) \right) \right) ( \nu ) = \frac{ \frac{
\zeta ( \nu ) }{ \xi ( \nu ) } }{ \frac{1}{ \xi ( \nu ) } } =
\zeta ( \nu )$, whenever $\nu \in {\bf N}$, while $\frac{1}{ \xi }
\in {\Gamma }_{ {\Phi }_{P} \left( \xi \right) }$, since
$\sum\limits_{ \nu = 0 }^{ \infty } \frac{1}{ \xi ( \nu ) } \left(
{\Phi }_{P} \left( \xi \right) \right) \left( \{ \nu \} \right) =
\sum\limits_{ \nu = 0 }^{ \infty } \frac{1}{ \xi ( \nu ) } \xi (
\nu ) P \left( \{ \nu \} \right) = 1$ and $\frac{1}{ \xi ( \nu ) }
> 0$ for every $n \in {\bf N}$. So ${\Psi }_{P , \xi } : {\Gamma
}_{P} \rightarrow {\Gamma }_{ {\Phi }_{P} \left( \xi \right) }$ is
an isometry. Moreover, if $\zeta \in {\Gamma }_{P}$, then $\left(
{\Phi }_{ {\Phi }_{P} \left( \xi \right) } \left( {\Psi }_{P , \xi
} \left( \zeta \right) \right) \right) \left( \{ \nu \} \right) =
\left( {\Psi }_{P , \xi } \left( \zeta \right) \right) ( \nu )
\left( {\Phi }_{P} \left( \xi \right) \right) \left( \{ \nu \}
\right) = \frac{ \zeta ( \nu ) }{ \xi ( \nu ) } \xi ( \nu ) P
\left( \{ \nu \} \right) = \zeta ( \nu ) P \left( \{ \nu \}
\right) = \left( {\Phi }_{P} \left( \zeta \right) \right) \left(
\{ \nu \} \right) $, whenever $\nu \in {\bf N}$, hence ${\Phi }_{
{\Phi }_{P} \left( \xi \right) } \left( {\Psi }_{P , \xi } \left(
\zeta \right) \right) = {\Phi }_{P} \left( \zeta \right) $. \hfill
$\bigtriangleup $

\section{The proof of 1.7}

\subsection{The proof of (i) in 1.7}

It is well-known that ${\bf C}^{*}$ constitutes a commutative
Polish group under multiplication and if $d(x,y) = \vert x - y
\vert + \left\vert \frac{1}{x} - \frac{1}{y} \right\vert $,
whenever $x$ and $y$ are in ${\bf C}^{*}$, then $d$ constitutes a
complete compatible metric on ${\bf C}^{*}$. Given any ${\bf g}
\in {\bf H}$ and any ${\bf h} \in {\bf H}$, we set $\rho \left(
{\bf g} , {\bf h} \right) = \sup\limits_{ n \in {\bf N} } d \left(
{\bf g}(n), {\bf h}(n) \right) $ and it is not difficult to verify
that $\rho $ constitutes a metric on ${\bf H}$. So let $\left(
{\bf h}_{k} \right) _{k \in {\bf N}}$ be any Cauchy sequence in
$\left( {\bf H}, \rho \right) $ and let $\epsilon > 0$. Then there
exists $K \in {\bf N}$ such that for any integer $k \geq K$ and
for any integer $l \geq K$, we have $\left\vert {\bf h}_{k}(n) -
{\bf h}_{l}(n) \right\vert \leq d \left( {\bf h}_{k}(n) , {\bf
h}_{l}(n) \right) \leq \rho \left( {\bf h}_{k} , {\bf h}_{l}
\right) < \frac{ \epsilon }{2} $, whenever $n \in {\bf N}$. So for
any $n \in {\bf N}$, $\left( {\bf h}_{k}(n) \right) _{k \in {\bf
N}}$ constitutes a Cauchy sequence in $\left( {\bf C}^{*} , d
\right) $ and consequently it has a limit, say ${\bf h}(n) =
\lim\limits_{k \rightarrow \infty }{\bf h}_{k}(n)$. Moreover,
since $\lim\limits_{n \rightarrow \infty } {\bf h}_{K}(n) = 1$,
there exists $N \in {\bf N}$ such that for any integer $n \geq N$,
we have $\left\vert {\bf h}_{K}(n) - 1 \right\vert < \frac{
\epsilon }{2}$ and hence $\left\vert {\bf h}(n) - 1 \right\vert =
\lim\limits_{l \rightarrow \infty } \left\vert {\bf h}_{l}(n) - 1
\right\vert \leq \sup\limits_{l \geq K} \left( \left\vert {\bf
h}_{l}(n) - {\bf h}_{K}(n) \right\vert + \left\vert {\bf h}_{K}(n)
- 1 \right\vert \right) \leq \sup\limits_{l \geq K} \left\vert
{\bf h}_{l}(n) - {\bf h}_{K}(n) \right\vert + \left\vert {\bf
h}_{K}(n) - 1 \right\vert < \epsilon $, while for any integer $k
\geq K$ and for any $n \in {\bf N}$, we have $d \left( {\bf
h}_{k}(n) , {\bf h}(n) \right) = \lim\limits_{l \rightarrow \infty
} d \left( {\bf h}_{k}(n) , {\bf h}_{l}(n) \right) \leq \frac{
\epsilon }{2}$, which implies that ${\bf h} \in {\bf H}$, hence
$\rho \left( {\bf h}_{k} , {\bf h} \right) = \sup\limits_{n \in
{\bf N}} d \left( {\bf h}_{k}(n) , {\bf h}(n) \right) \leq \frac{
\epsilon }{2} < \epsilon $ and consequently ${\bf h}_{k}
\rightarrow {\bf h}$ in $\left( {\bf H} , \rho \right) $ as $k
\rightarrow \infty $, which implies that $\rho $ constitutes a
complete metric on ${\bf H}$. If ${\bf f}$, ${\bf g}$ and ${\bf
h}$ are any elements of ${\bf H}$, then it is not difficult to
prove that $\rho \left( \frac{1}{\bf f} , \frac{1}{\bf g} \right)
= \rho \left( {\bf f} , {\bf g} \right) $, which implies that
inversion is continuous, and $$\rho \left( {\bf f}{\bf h} , {\bf
g}{\bf h} \right) \leq \max \left\{ \sup\limits_{n \in {\bf N}}
\left\vert {\bf h}(n) \right\vert , \sup\limits_{n \in {\bf N}}
\frac{1}{ \left\vert {\bf h}(n) \right\vert } \right\} \rho \left(
{\bf f} , {\bf g} \right) ,$$ since $\lim\limits_{n \rightarrow
\infty } {\bf h}(n) = 1$ and the same holds true for $\frac{1}{\bf
h}$. So let ${\bf f}_{k} \rightarrow {\bf f}$ in $\left( {\bf H} ,
\rho \right) $ as $k \rightarrow \infty $ and let ${\bf g}_{k}
\rightarrow {\bf g}$ in $\left( {\bf H} , \rho \right) $ as $k
\rightarrow \infty $, while $\epsilon > 0$. Then there exists $K
\in {\bf N}$ such that for any integer $k \geq K$, we have
$$\rho \left( {\bf g}_{k} , {\bf g} \right) < \frac{ \epsilon }{ 2 \max
\left\{ 1 + \sup\limits_{n \in {\bf N}} \left\vert {\bf f}(n)
\right\vert , 1 + \sup\limits_{n \in {\bf N}} \frac{1}{ \left\vert
{\bf f}(n) \right\vert } \right\} }$$ and $$\rho \left( {\bf
f}_{k} , {\bf f} \right) < \min \left\{ \frac{ \epsilon }{ 2 \max
\left\{ \sup\limits_{n \in {\bf N}} \left\vert {\bf g}(n)
\right\vert , \sup\limits_{n \in {\bf N}} \frac{1}{ \left\vert
g(n) \right\vert } \right\} } , 1 \right\} ,$$ so $$\sup\limits_{n
\in {\bf N}} \left\vert {\bf f}_{k}(n) \right\vert \leq
\sup\limits_{n \in {\bf N}} \left\vert {\bf f}_{k}(n) - {\bf f}(n)
\right\vert + \sup\limits_{n \in {\bf N}} \left\vert {\bf f}(n)
\right\vert < 1 + \sup\limits_{n \in {\bf N}} \left\vert {\bf
f}(n) \right\vert $$ and $$\sup\limits_{n \in {\bf N}}\frac{1}{
\left\vert {\bf f}_{k}(n) \right\vert } \leq \sup\limits_{n \in
{\bf N}} \left\vert \frac{1}{{\bf f}_{k}(n)} - \frac{1}{{\bf
f}(n)} \right\vert + \sup\limits_{n \in {\bf N}}\frac{1}{
\left\vert {\bf f}(n) \right\vert } < 1 + \sup\limits_{n \in {\bf
N}}\frac{1}{ \left\vert {\bf f}(n) \right\vert },$$ so
\begin{enumerate}
\item[ ]
$\rho \left( {\bf f}_{k}{\bf g}_{k} , {\bf f}{\bf g} \right) \leq
\rho \left( {\bf f}_{k}{\bf g}_{k} , {\bf f}_{k}{\bf g} \right) +
\rho \left( {\bf f}_{k}{\bf g} , {\bf f}{\bf g} \right) $
\item[ ]
$\leq \max \left\{ \sup\limits_{n \in {\bf N}} \left\vert {\bf
f}_{k}(n) \right\vert , \sup\limits_{n \in {\bf N}} \frac{1}{
\left\vert {\bf f}_{k}(n) \right\vert } \right\} \rho \left( {\bf
g}_{k} , {\bf g} \right) $
\item[ ]
$+ \max \left\{ \sup\limits_{n \in {\bf N}} \left\vert {\bf g}(n)
\right\vert , \sup\limits_{n \in {\bf N}}\frac{1}{ \left\vert {\bf
g}(n) \right\vert } \right\} \rho \left( {\bf f}_{k} , {\bf f}
\right) $
\item[ ]
$< \epsilon $.
\end{enumerate}
So ${\bf H}$ constitutes a topological group whose topology is
given by the complete metric $\rho $. What is left to show is that
$\left( {\bf H} , \rho \right) $ is separable. Indeed, it is not
difficult to verify that $${\mathcal{C}} = \left\{ {\bf g} \in
\left( \left( {\bf Q} + i {\bf Q} \right) \setminus \{ 0 \}
\right) ^{\bf N} : \exists m \forall n \geq m \left( {\bf g}(n) =
1 \right) \right\} $$ constitutes a countable dense subset of
$\left( {\bf H} , \rho \right) $. \hfill $\bigtriangleup $

\subsection{The proof of (ii) in 1.7}

If ${\bf h} \in {\bf H}$ and ${\bf x} \in {\ell }^{1} \left( {\bf
C}^{*} \right) $, then
\begin{enumerate}
\item[ ]
$\left\Vert {\bf h} \cdot {\bf x} \right\Vert _{1} =
\sum\limits_{n=0}^{ \infty } \left\vert {\bf h}(n){\bf x}(n)
\right\vert $
\item[ ]
$\leq \left( \sup\limits_{n \in {\bf N}} \left\vert {\bf h}(n)
\right\vert \right) \sum\limits_{n=0}^{ \infty } \left\vert {\bf
x}(n) \right\vert $
\item[ ]
$= \left( \sup\limits_{n \in {\bf N}} \left\vert {\bf h}(n)
\right\vert \right) \left\Vert {\bf x} \right\Vert _{1}$
\end{enumerate}
and since ${\bf h} \in {\bf H}$, the fact that ${\bf x} \in {\ell
}^{1} \left( {\bf C}^{*} \right) $, implies that ${\bf h} \cdot
{\bf x} \in {\ell }^{1} \left( {\bf C}^{*} \right) $ and it is not
difficult to verify that $\left( {\bf h} , {\bf x} \right) \mapsto
{\bf h} \cdot {\bf x}$ is a group action. So let ${\bf h}_{k}
\rightarrow {\bf h}$ in ${\bf H}$ as $k \rightarrow \infty $ and
${\bf x}_{k} \rightarrow {\bf x}$ in ${\ell }^{1} \left( {\bf
C}^{*} \right) $ as $k \rightarrow \infty $. Then
\begin{enumerate}
\item[ ]
$\left\Vert {\bf h}_{k} \cdot {\bf x}_{k} - {\bf h} \cdot {\bf x}
\right\Vert _{1} \leq \left\Vert {\bf h}_{k} \cdot \left( {\bf
x}_{k} - {\bf x} \right) \right\Vert _{1} + \left\Vert \left( {\bf
h}_{k} - {\bf h} \right) \cdot {\bf x} \right\Vert _{1}$
\item[ ]
$\leq \left( \sup\limits_{n \in {\bf N}} \left\vert {\bf h}_{k}(n)
\right\vert \right) \left\Vert {\bf x}_{k} - {\bf x} \right\Vert
_{1} + \sup\limits_{n \in {\bf N}} \left\vert {\bf h}_{k}(n) -
{\bf h}(n) \right\vert \Vert {\bf x} \Vert _{1}$
\item[ ]
$\leq \left( \sup\limits_{n \in {\bf N}} \left\vert {\bf h}_{k}(n)
\right\vert \right) \left\Vert {\bf x}_{k} - {\bf x} \right\Vert
_{1} + \rho \left( {\bf h}_{k} , {\bf h} \right) \Vert {\bf x}
\Vert _{1}$
\item[ ]
$\leq \left( \sup\limits_{n \in {\bf N}} \left\vert {\bf h}_{k}(n)
- {\bf h}(n) \right\vert + \sup\limits_{n \in {\bf N}} \left\vert
{\bf h}(n) \right\vert \right) \left\Vert {\bf x}_{k} - {\bf x}
\right\Vert _{1} + \rho \left( {\bf h}_{k} , {\bf h} \right) \Vert
{\bf x} \Vert _{1}$
\item[ ]
$\leq \left( \rho \left( {\bf h}_{k} , {\bf h} \right) +
\sup\limits_{n \in {\bf N}} \left\vert {\bf h}(n) \right\vert
\right) \left\Vert {\bf x}_{k} - {\bf x} \right\Vert _{1} + \rho
\left( {\bf h}_{k} , {\bf h} \right) \Vert {\bf x} \Vert _{1}
\rightarrow 0$
\end{enumerate}
as $k \rightarrow \infty $. \hfill $\bigtriangleup $

\subsection{The proof of (iii) in 1.7}

It is enough to notice that if ${\bf y} \in {\ell }^{1} \left(
{\bf C}^{*} \right) $ and $N \in {\bf N}$, while
\[
{\bf h}_{N}(n) = \left\{
\begin{array}{ll}
\frac{ {\bf y}(n) }{ {\bf x}(n) } & \mbox{if $n \in \{ 0, ..., N
\} $}
\\ \\
1 & \mbox{if $n \in {\bf N} \setminus \{ 0, ..., N \} $}
\end{array}
\right.
\]
then ${\bf h}_{N} \in {\bf H}$ and $\left\Vert {\bf h}_{N} \cdot
{\bf x} - {\bf y} \right\Vert _{1} = \sum\limits_{n > N}
\left\vert {\bf x}(n) - {\bf y}(n) \right\vert \rightarrow 0$ as
$N \rightarrow \infty $. \hfill $\bigtriangleup $

\subsection{The proof of (iv) in 1.7}

If ${\bf y} \in {\bf H} \cdot {\bf x}$, then it is not difficult
to verify that $\lim\limits_{n \rightarrow \infty } \frac{ {\bf
y}(n) }{ {\bf x}(n) } = 1$ and consequently there exists $m \in
{\bf N}$ such that for any integer $n \geq m$, we have $\left\vert
\frac{ {\bf y}(n) }{ {\bf x}(n) } \right\vert \leq \frac{3}{2}$.
So ${\bf H} \cdot {\bf x} \subseteq {\mathcal{M}}$, where
$${\mathcal{M}} = \bigcup\limits_{m \in {\bf N}} \bigcap\limits_{n
\geq m} \left\{ {\bf y} \in {\ell }^{1} \left( {\bf C}^{*} \right)
: \left\vert \frac{ {\bf y}(n) }{ {\bf x}(n) } \right\vert \leq
\frac{3}{2} \right\} $$ is easily verified to be $F_{\sigma }$. So
it is enough to prove that ${\ell }^{1} \left( {\bf C}^{*} \right)
\setminus {\mathcal{M}}$ is dense in ${\ell }^{1} \left( {\bf
C}^{*} \right) $. Indeed, if ${\bf z} \in {\ell }^{1} \left( {\bf
C}^{*} \right) $ and $N \in {\bf N}$, while
\[
{\bf z}_{N}(n) = \left\{
\begin{array}{ll}
{\bf z}(n) & \mbox{if $n \in \{ 0, ..., N \} $}
\\ \\
2{\bf x}(n) & \mbox{if $n \in {\bf N} \setminus \{ 0, ..., N \} $}
\end{array}
\right.
\]
then it is enough to notice that ${\bf z}_{N} \in {\ell }^{1}
\left( {\bf C}^{*} \right) \setminus {\mathcal{M}}$ and
$$\left\Vert {\bf z}_{N} - {\bf z} \right\Vert _{1} = \sum\limits_{n >
N} \left\vert 2{\bf x}(n) - {\bf z}(n) \right\vert \rightarrow 0$$
as $N \rightarrow \infty $. \hfill $\bigtriangleup $

\section{The proof of 1.9}

\subsection{The proof of (i)}

It is well-known that $\left( 0 , \infty \right) $ constitutes a
commutative Polish group under multiplication and if $d(x,y) =
\vert x - y \vert + \left\vert \frac{1}{x} - \frac{1}{y}
\right\vert $, whenever $x$ and $y$ are in $\left( 0 , \infty
\right) $, then $d$ constitutes a complete compatible metric on
$\left( 0 , \infty \right) $. Given any $f \in {\bf F}$ and any $g
\in {\bf F}$, we set $\rho (f,g) = \sup\limits_{x>0}d \left( f(x),
g(x) \right) $ and it is not difficult to verify that $\rho $
constitutes a metric on ${\bf F}$. So let $\left( f_{k} \right)
_{k \in {\bf N}}$ be any Cauchy sequence in $\left( {\bf F}, \rho
\right) $ and let $\epsilon > 0$. Then there exists $K \in {\bf
N}$ such that for any integer $k \geq K$ and for any integer $l
\geq K$, we have $\left\vert f_{k}(x) - f_{l}(x) \right\vert \leq
d \left( f_{k}(x) , f_{l}(x) \right) \leq \rho \left( f_{k} ,
f_{l} \right) < \frac{ \epsilon }{2} $, whenever $x>0$. So for any
$x>0$, $\left( f_{k}(x) \right) _{k \in {\bf N}}$ constitutes a
Cauchy sequence in $\left( \left( 0, \infty \right) , d \right) $
and consequently it has a limit, say $f(x) = \lim\limits_{k
\rightarrow \infty }f_{k}(x)$. Moreover, since $\lim\limits_{x
\rightarrow 0} f_{K}(x) = 1$, there exists $\delta
>0$ such that for any $x \in (0 , \delta )$, we have $\left\vert
f_{K}(x) - 1 \right\vert < \frac{ \epsilon }{2}$ and hence
$\left\vert f(x) - 1 \right\vert = \lim\limits_{l \rightarrow
\infty } \left\vert f_{l}(x) - 1 \right\vert \leq \sup\limits_{l
\geq K} \left( \left\vert f_{l}(x) - f_{K}(x) \right\vert +
\left\vert f_{K}(x) - 1 \right\vert \right) \leq \sup\limits_{l
\geq K} \left\vert f_{l}(x) - f_{K}(x) \right\vert + \left\vert
f_{K}(x) - 1 \right\vert < \epsilon $, and, in addition, since
$\lim\limits_{x \rightarrow \infty } f_{K}(x) = 1$, there exists
$M>0$ such that for any $x \geq M$, we have $\left\vert f_{K}(x) -
1 \right\vert < \frac{ \epsilon }{2}$ and hence $\left\vert f(x) -
1 \right\vert = \lim\limits_{l \rightarrow \infty } \left\vert
f_{l}(x) - 1 \right\vert \leq \sup\limits_{l \geq K} \left(
\left\vert f_{l}(x) - f_{K}(x) \right\vert + \left\vert f_{K}(x) -
1 \right\vert \right) \leq \sup\limits_{l \geq K} \left\vert
f_{l}(x) - f_{K}(x) \right\vert + \left\vert f_{K}(x) - 1
\right\vert < \epsilon $, while for any integer $k \geq K$ and for
any $x>0$, we have $d \left( f_{k}(x) , f(x) \right) =
\lim\limits_{l \rightarrow \infty } d \left( f_{k}(x) , f_{l}(x)
\right) \leq \frac{ \epsilon }{2}$, which, given [4], implies that
$f \in {\bf F}$, hence $\rho \left( f_{k} , f \right) =
\sup\limits_{x>0} d \left( f_{k}(x) , f(x) \right) \leq \frac{
\epsilon }{2} < \epsilon $ and consequently $f_{k} \rightarrow f$
in $\left( {\bf F} , \rho \right) $ as $k \rightarrow \infty $,
which implies that $\rho $ constitutes a complete metric on ${\bf
F}$. If $f$, $g$ and $h$ are any elements of ${\bf F}$, then it is
not difficult to prove that $\rho \left( \frac{1}{f} , \frac{1}{g}
\right) = \rho (f,g)$, which implies that inversion is continuous,
and $$\rho (fh,gh) \leq \max \left\{ \sup\limits_{x>0}h(x) ,
\sup\limits_{x>0} \frac{1}{h(x)} \right\} \rho (f,g),$$ since
$\lim\limits_{x \rightarrow 0} h(x) = \lim\limits_{x \rightarrow
\infty } h(x) = 1$ and $h$ being continuous on any compact
interval it is also bounded, while the same hold true for
$\frac{1}{h}$. So let $f_{k} \rightarrow f$ in $\left( {\bf F} ,
\rho \right) $ as $k \rightarrow \infty $ and let $g_{k}
\rightarrow g$ in $\left( {\bf F} , \rho \right) $ as $k
\rightarrow \infty $, while $\epsilon > 0$. Then there exists $K
\in {\bf N}$ such that for any integer $k \geq K$, we have
$$\rho \left( g_{k} , g \right) < \frac{ \epsilon }{ 2 \max
\left\{ 1 + \sup\limits_{x>0}f(x) , 1 + \sup\limits_{x>0}
\frac{1}{f(x)} \right\} }$$ and $$\rho \left( f_{k} , f \right) <
\min \left\{ \frac{ \epsilon }{ 2 \max \left\{
\sup\limits_{x>0}g(x) , \sup\limits_{x>0} \frac{1}{g(x)} \right\}
} , 1 \right\} ,$$ so $$\sup\limits_{x>0}f_{k}(x) \leq
\sup\limits_{x>0} \left\vert f_{k}(x) - f(x) \right\vert +
\sup\limits_{x>0}f(x) < 1 + \sup\limits_{x>0}f(x)$$ and
$$\sup\limits_{x>0}\frac{1}{f_{k}(x)} \leq \sup\limits_{x>0}
\left\vert \frac{1}{f_{k}(x)} - \frac{1}{f(x)} \right\vert +
\sup\limits_{x>0}\frac{1}{f(x)} < 1 +
\sup\limits_{x>0}\frac{1}{f(x)},$$ so
\begin{enumerate}
\item[ ]
$\rho \left( f_{k}g_{k} , fg \right) \leq \rho \left( f_{k}g_{k} ,
f_{k}g \right) + \rho \left( f_{k}g , fg \right) $
\item[ ]
$\leq \max \left\{ \sup\limits_{x>0}f_{k}(x) , \sup\limits_{x>0}
\frac{1}{f_{k}(x)} \right\} \rho \left( g_{k} , g \right) $
\item[ ]
$+ \max \left\{ \sup\limits_{x>0}g(x) ,
\sup\limits_{x>0}\frac{1}{g(x)} \right\} \rho \left( f_{k} , f
\right) $
\item[ ]
$< \epsilon $.
\end{enumerate}
So ${\bf F}$ constitutes a topological group whose topology is
given by the complete metric $\rho $ and what is left to show is
that $\left( {\bf F} , \rho \right) $ is separable. Given any
integer $N \geq 2$, we denote by ${\mathcal{C}}_{N}$ the set of
all $\phi $ with the property that there exists $P(X) \in {\bf
Q}[X]$ with positive values on $\left[ \frac{1}{N} , N \right] $
such that $\phi \left\vert \left[ \frac{1}{N} , N \right] \right.
= P \left\vert \left[ \frac{1}{N} , N \right] \right. $ and $\phi
\left\vert \left( \left( 0, \frac{1}{2N} \right] \cup \left[ N +
\frac{1}{N} , \infty \right) \right) = 1 \right. $, while
$\frac{1}{2N} \leq x \leq \frac{1}{N} \Rightarrow \phi (x) = N
\left( x - \frac{1}{N} \right) \left( 1 - P \left( \frac{1}{N}
\right) \right) + P \left( \frac{1}{N} \right) $ and $N \leq x
\leq N + \frac{1}{N} \Rightarrow \phi (x) = N(x-N) \left( 1 - P(N)
\right) + P(N) $, and let ${\mathcal{C}} = \bigcup\limits_{N \geq
2}{\mathcal{C}}_{N}$. It is not difficult to verify that
${\mathcal{C}}$ is countable and, given [4], that ${\mathcal{C}}$
is dense in $\left( {\bf F} , \rho \right) $. \hfill
$\bigtriangleup $

\subsection{The proof of (ii)}

If $f \in {\bf F}$ and $g \in L^{1}_{++} \left( (0, \infty ) ,
\lambda \right) $, then $$\left\Vert fg \right\Vert _{1} =
\int_{0}^{ \infty } f(x)g(x)dx \leq \left( \sup\limits_{x>0}f(x)
\right) \int_{0}^{ \infty }g(x)dx = \left( \sup\limits_{x>0}f(x)
\right) \left\Vert g \right\Vert _{1}$$ and since $f>0$, the fact
that $g>0$, $\lambda $-a.e. implies that $fg>0$, $\lambda $-a.e.,
so $fg \in L^{1}_{++} \left( (0, \infty ) , \lambda \right) $ and
it is not difficult to verify that $(f,g) \mapsto fg$ is a group
action. So let $f_{k} \rightarrow f$ in ${\bf F}$ as $k
\rightarrow \infty $ and $g_{k} \rightarrow g$ in $L^{1}_{++}
\left( (0, \infty ) , \lambda \right) $ as $k \rightarrow \infty
$. Then
\begin{enumerate}
\item[ ]
$\left\Vert f_{k}g_{k} - fg \right\Vert _{1} \leq \left\Vert f_{k}
\left( g_{k} - g \right) \right\Vert _{1} + \left\Vert \left(
f_{k} - f \right) g \right\Vert _{1}$
\item[ ]
$\leq \left( \sup\limits_{x>0} f_{k}(x) \right) \left\Vert g_{k} -
g \right\Vert _{1} + \sup\limits_{x>0} \left\vert f_{k}(x) - f(x)
\right\vert \Vert g \Vert _{1}$
\item[ ]
$\leq \left( \sup\limits_{x>0} f_{k}(x) \right) \left\Vert g_{k} -
g \right\Vert _{1} + \rho \left( f_{k} , f \right) \Vert g \Vert
_{1}$
\item[ ]
$\leq \left( \sup\limits_{x>0} \left\vert f_{k}(x) - f(x)
\right\vert + \sup\limits_{x>0}f(x) \right) \left\Vert g_{k} - g
\right\Vert _{1} + \rho \left( f_{k} , f \right) \Vert g \Vert
_{1}$
\item[ ]
$\leq \left( \rho \left( f_{k} , f \right) + \sup\limits_{x>0}f(x)
\right) \left\Vert g_{k} - g \right\Vert _{1} + \rho \left( f_{k}
, f \right) \Vert g \Vert _{1} \rightarrow 0$
\end{enumerate}
as $k \rightarrow \infty $. \hfill $\bigtriangleup $

\section{The proof of 1.12}

We denote by ${\rho }^{*}$ the complete metric on ${\bf G}^{*}$.
If ${\bf g}$, ${\bf h}$ are any points of ${\bf G}^{*}$ such that
$f_{\bf g} = f_{\bf h}$, then for any $n \in {\bf N}^{*}$, we have
that ${\bf g}(n) = f_{\bf g}(n) = f_{\bf h}(n) = {\bf h}(n)$, so
${\bf g} = {\bf h}$ and consequently we have an injection. So if
${\bf g}_{k} \rightarrow {\bf g}$ in ${\bf G}^{*}$ as $k
\rightarrow \infty $, then given $n \in {\bf N}^{*}$ and $x \in
[n,n+1]$, we have that $\left\vert f_{{\bf g}_{k}}(x) - f_{\bf
g}(x) \right\vert \leq \left\vert {\bf g}_{k}(n+1) - {\bf g}(n+1)
\right\vert \cdot \vert x-n \vert + \left\vert {\bf g}_{k}(n) -
{\bf g}(n) \right\vert \cdot \vert x - (n+1) \vert \leq 2 {\rho
}^{*} \left( {\bf g}_{k} , {\bf g} \right) $, while given $x \in
\left[ \frac{1}{2} , 1 \right] $, we have that $\left\vert f_{{\bf
g}_{k}}(x) - f_{\bf g}(x) \right\vert \leq \left\vert {\bf
g}_{k}(1) - {\bf g}(1) \right\vert \cdot \vert 2x-1 \vert \leq
{\rho }^{*} \left( {\bf g}_{k} , {\bf g} \right) $, hence
$\sup\limits_{x>0} \left\vert f_{{\bf g}_{k}}(x) - f_{\bf g}(x)
\right\vert \leq 2 {\rho }^{*} \left( {\bf g}_{k} , {\bf g}
\right) \rightarrow 0$ as $k \rightarrow \infty $ and consequently
$\rho \left( f_{{\bf g}_{k}} , f_{\bf g} \right) \rightarrow 0$ as
$k \rightarrow \infty $. So the injection in question is
continuous, but if ${\bf g}$, ${\bf h}$ in ${\bf G}^{*}$ are such
that ${\bf g}(1) \neq 1$, ${\bf h}(1) \neq 1$, then for any $x \in
\left( 0 , \frac{1}{2} \right) $, we have that $f_{\bf g}(x)
f_{\bf h}(x) - f_{{\bf g}{\bf h}}(x) = 2  (x-1) \cdot (2x-1) \cdot
\left( {\bf g}(1)-1 \right) \cdot \left( {\bf h}(1)-1 \right) $
and consequently $f_{\bf g}(x)f_{\bf h}(x) \neq f_{{\bf g}{\bf
h}}(x)$. \hfill $\bigtriangleup $

\end{document}